\documentclass[twocolumn]{autart}    

\usepackage{graphicx}          
\usepackage{amsmath}
\usepackage{amssymb}
\usepackage{color}

\usepackage{savesym} \savesymbol{AND}

\usepackage{algorithmic}
\usepackage{algorithm}

\newtheorem{theorem}{Theorem}
\newtheorem{lemma}[theorem]{Lemma}
\newtheorem{assumption}[theorem]{Assumption}
\newtheorem{remark}[theorem]{Remark}

\allowdisplaybreaks
 
\newcommand{\R}{\mathbb{R}}
\newcommand{\N}{\mathbb{N}}
\newcommand{\A}{\mathcal{A}}
\newcommand{\B}{\mathcal{B}}
\newcommand{\D}{\mathcal{D}}

\begin{document}

\begin{frontmatter}

\title{Predictive feedback boundary control of semilinear and quasilinear $2\times 2$ hyperbolic PDE-ODE systems  \thanksref{footnoteinfo}} 
\thanks[footnoteinfo]{This paper was not presented at any IFAC 
meeting. Corresponding author Timm Strecker.}

\author[Melb]{Timm Strecker}\ead{timm.strecker@unimelb.edu.au},    
\author[NTNU]{Ole Morten Aamo}\ead{ole.morten.aamo@ntnu.no},    
\author[Melb]{Michael Cantoni}\ead{cantoni@unimelb.edu.au}

\address[Melb]{Department of Electrical and Electronic Engineering, The University of
Melbourne, Parkville 3010, Australia}  
\address[NTNU]{Department of Engineering Cybernetics, Norwegian University
of Science and Technology (NTNU), Trondheim N-7491, Norway}  


\begin{abstract}                          
We present a control design for semilinear and quasilinear $2\times 2$ hyperbolic partial differential equations with the  control input  at one boundary and a nonlinear ordinary differential equation coupled to the other. The controller can be designed to asymptotically stabilize the system at an equilibrium or relative to a reference signal.  Two related but different controllers for semilinear and general quasilinear systems are presented and the additional challenges in quasilinear systems are discussed. Moreover, we  present an observer that estimates the distributed PDE state and the unmeasured ODE state from measurements at the actuated boundary only, which can be used to  also solve the output feedback control problem.
\end{abstract}

\end{frontmatter}

\section{Introduction}
\label{sec:introduction}

In this paper, we consider systems consisting of a 1-d  hyperbolic partial differential equation (PDE) coupled with an ordinary differential equation (ODE) at the uncontrolled boundary, as given by
%
%
%
\begin{align}
w_t(x,t) &= \Lambda(x,w(x,t))\,w_x(x,t) + F(x,w(x,t)), \label{w_t}\\
\dot{X}(t) &= f^0(X(t),v(0,t),t), \label{X_t}\\
u(0,t) &= g^0(X(t),v(0,t),t), \label{uBC}\\
v(1,t) & = U(t), \label{vBC}\\
w(\cdot,0)&=w_0=\left( \begin{matrix}u_0&v_0 \end{matrix} \right)^T,  \label{uvIC} \\
X(0) &= X_0. \label{XIC}
\end{align}
Here $w(x,t)=\left(\begin{matrix} u(x,t)&v(x,t) \end{matrix} \right)^T$, $X(t)\in\R^{n}$,
subscripts denote partial derivatives,  $x\in[0,1]$, $t\geq 0$, $w_0$ and $X_0$ are initial conditions,  and $U$ is the control input.
 The nonlinear functions $\Lambda$ and $F$, which model the propagation speeds and source terms, respectively, have the structure 
\begin{align}
\Lambda(x,z) &=\left( \begin{matrix}-\lambda^u(x,z) & 0\\ 0& \lambda^v(x,z)    \end{matrix}  \right), & \lambda^u,\,\lambda^v&>0, \\
F(x,z) &=\left( \begin{matrix} f^u(x,z)& f^v(x,z)    \end{matrix}\right)^T.
\end{align}
The nonlinear functions $f^0$ and $g^0$ define the coupling to the ODE at the uncontrolled boundary.
More precise assumptions on the system  data are given in the subsequent sections.

The system is called \emph{semilinear} in the special case where $\Lambda$ is independent of the state $w(x,t)$, whereas a  \emph{quasilinear} system is the general case when $\Lambda$ is a function of the state $w(x,t)$. The difference between semilinear and quasilinear systems is significant. In particular, the characteristic lines of quasilinear systems depend on the state and the control input. If characteristic lines ``collide'', the solution ceases to exist due to blow-up of the gradient of that state, even if the state itself remains bounded \cite{bressan2000hyperbolic}. Therefore, the control inputs need to be designed not only to regulate the state, but also to prevent a collision of characteristic lines.  By contrast, the characteristic lines of semilinear (including linear) systems are known a priori and they do not collide, which simplifies the control design. 
The dependence of $\Lambda$ on the state also has implication on the state space. For semilinear systems, so called broad solutions\footnote{Broad solutions are a type of weak solution that is defined as the solution of the integral equations that are obtained by integrating (\ref{w_t}) along its characteristic lines.} \cite[Chapter 3]{bressan2000hyperbolic} can be defined over $L^{\infty}$. That is, the solutions can be discontinuous and the control inputs can be chosen freely in $L^{\infty}$. In quasilinear systems, at least Lipschitz-continuity of the state is required in order to define broad solutions. As a consequence, the control inputs must be compatible with the current state at all times, and there is a limit on how fast the system can be steered to the desired state. 

Hyperbolic systems similar to (\ref{w_t}) have received significant attention. They model several 1-d transport phenomena such as gas or fluid flow through pipelines, open-channel flow, traffic flow and electrical transmission lines \cite{bastin2016book}. Control design approaches for these systems  include dissipative boundary conditions \cite{greenberg1984effect} and control Lyapunov functions \cite{coron2007strict}. However, the use of such static boundary feedback   works only if the source terms are small \cite{bastin2016book}. An alternative approach for linear systems that steers the system more actively to the desired state is backstepping, which now exists for a variety of system classes, see e.g.~\cite{krstic2008backstepping,vazquez2011backstepping,coron2013local}. Motivated by an application to mechanical vibrations in drill strings \cite{di2015distributed}, a backstepping controller  for PDE-ODE systems  similar to those considered here is developed in   \cite{di2018PDEODE}, although only  for linear systems. A different approach has been pursued in \cite{strecker2017output,strecker2019semilineargeneralheterodirectional} for semilinear systems and in \cite{strecker2019quasilinearfirstorder,strecker2021quasilinear2x2} for quasilinear systems, which achieves  performance equivalent to that of backstepping controllers, but  for nonlinear systems.   In \cite{strecker2021quasilinear2x2}, robustness to parameter uncertainty and measurement and actuation errors is also investigated.

The contributions of this paper extend  \cite{strecker2021quasilinear2x2} to  hyperbolic systems that are coupled with a nonlinear ODE at the uncontrolled boundary, as given by (\ref{w_t})-(\ref{XIC}). Moreover, we present an observer that estimates the distributed PDE state $w$ and the ODE state $X$ from boundary measurements of $u$ at the actuated boundary only.

The paper is organized as follows. First, the state feedback control design for semilinear systems is presented in Section \ref{sec: semilinear}. The  control design for quasilinear systems is given in Section \ref{sec: quasilinear}. The observer  is presented in Section \ref{sec:observer}. In Section \ref{sec:simulation}, a simulation study is performed to demonstrate the effectiveness of the proposed controller. Section \ref{sec:conclusion} contains some concluding remarks. 

\section{Semilinear systems} \label{sec: semilinear}
The control design proposed here for semilinear systems builds on ideas first presented in \cite{strecker2017output}. The approach is based on the following two observations:
\begin{itemize}
\item Due to the hyperbolic nature of Equation (\ref{w_t}), the effect of the control input entering in boundary condition (\ref{vBC}) propagates through the domain with finite speed $\lambda^v$. Therefore, the state $w$ in the interior of the domain and the state $X$ of the ODE (\ref{X_t}), are only affected after a certain amount of time has passed.
\item  The system is easier to control by  treating $v(0,\cdot)$, the  value  of $v$ at the uncontrolled boundary,  as a  ``virtual'' control input $U^*$. This is because $v(0,\cdot)$ enters  both the boundary condition for $u$, as given in (\ref{uBC}), and in the ODE (\ref{X_t}), whereas the coupling between $U$ and the states $u$ and $X$ is more indirect and  affected by the delay.  In a second step, the actual inputs $U$ can be constructed such that $v(0,\cdot)$ becomes equal to $U^*$. 
\end{itemize}
For semilinear systems, we present a continuous-time state feedback control law. At each time $t$ the controller maps the  state $(w(\cdot,t),\,X(t))$ into the control input $U(t)$. This state feedback controller can be combined with the observer from Section \ref{sec:observer} to obtain an output feedback controller.

\subsection{Assumptions}
For semilinear systems considered in this section, we make the following assumptions.
\begin{assumption}\label{assumption semilinear}
The speeds in $\Lambda$ are assumed to be independent of the state $w$. The speeds are bounded from below by
\begin{equation}
\sup_{x\in[0,1]} \min\{ (\lambda^u(x))^{-1},\,(\lambda^v(x))^{-1} \} \leq l_{\Lambda^{-1}}.
\end{equation}
 The nonlinear functions $F$, $f^0$, $g^0$ are globally Lipschitz-continuous in the state arguments, i.e.,
\begin{align}
 \operatorname*{ess\,sup}_{x\in[0,1],w\in\R^2}  \|\partial_w F(x, w) \|_{\infty} &=l_F,  \label{boundF_w}\\
\operatorname*{ess\,sup}_{X\in\R^n, v \in\R,t\geq0}  |\partial_X f^0(X,v,t) | &= l_{f^0_X},  \\
\operatorname*{ess\,sup}_{X\in\R^n, v \in\R,t\geq0}  |\partial_v f^0(X,v,t) | &= l_{f^0_v},  \\
\operatorname*{ess\,sup}_{X\in\R^n, v \in\R,t\geq0}  |\partial_X g^0(X,v,t) | &= l_{g^0_X},  \\
\operatorname*{ess\,sup}_{X\in\R^n, v \in\R,t\geq0}  |\partial_v g^0(X,v,t) | &= l_{g^0_v},  \label{bound g0_v}
\end{align}
where $\partial$ denotes partial derivatives\footnote{By Rademacher's theorem, the partial derivatives of Lipschitz-continuous functions exist almost everywhere \cite[Theorem 2.8]{bressan2000hyperbolic}.} and $l_F$ through $ l_{g^0_v}$ are the finite Lipschitz constants. We assume further that 
\begin{align}
F(x,0)&=0 & \text{ for all } &x\in[0,1],\label{F(0)=0}  \\
 f^0(0,0,t)&=g^0(0,0,t) = 0 & \text{ for all } &t\geq 0,  \label{f0(0)=0}
\end{align}
which ensures that the origin is an equilibrium. Finally, assume there exists a controller $K(X(t),t)$ satisfying 
\begin{align}
\operatorname*{ess\,sup}_{X\in\R^n, t\geq0}  |\partial_X K(X,t) | &= l_{K}, \label{bound K_X}\\
 K(0,t)&= 0 & \text{ for all } &t\geq 0, \label{K(0)=0}
\end{align}
 such that the closed-loop system
\begin{equation}
\dot{X}(t) = f^0(X(t),K(X(t),t),t) \label{closed loop stable semilinear}
\end{equation}
has a  globally asymptotically stable equilibrium at the origin.
\end{assumption}

\subsection{Characteristic lines and determinate sets}
\begin{figure}[htbp!]\centering
\includegraphics[width=.6\columnwidth]{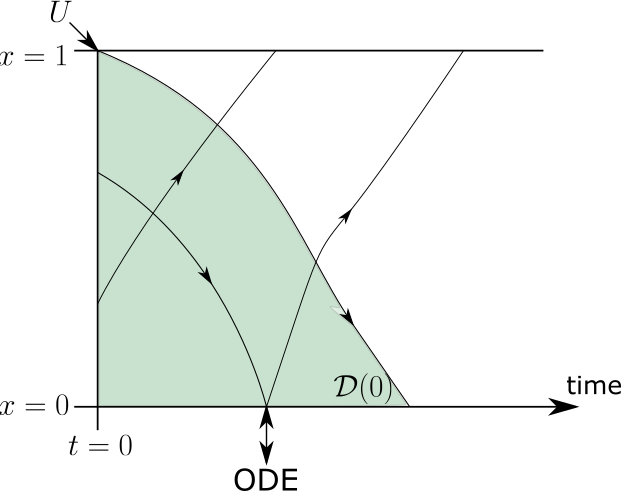}
\caption{Characteristic lines of system (\ref{w_t})-(\ref{XIC}), with the input $U$ entering at $x = 1$ and the coupling with the ODE at $x = 0$. The determinate set $\D(0)$ is shaded in green.}
\label{fig: characteristic lines}
\end{figure}
The characteristic lines of the system are sketched in Figure \ref{fig: characteristic lines}.
We parametrize  the characteristic lines of (\ref{w_t}) via
\begin{align}
\tau^u(t;x) &= t+\int_0^x\frac{1}{\lambda^u(\xi,w(\xi,\tau^u(t;\xi))}d\xi \label{tau^u},\\
\tau^v(t;x) &= t+\int_x^1\frac{1}{\lambda^v(\xi,w(\xi,\tau^v(t;\xi))}d\xi, \label{tau^v} 
\end{align}
where $\tau^u$ and $\tau^v$ are the times at which the characteristic lines that start at time $t$ and the spatial boundary $x=0$ and $x=1$, respectively, reach the location $x\in[0,1]$. 
The corresponding closed and half-open determinate sets are defined as
\begin{align}
 \D(t) &= \left\{(x,s):\,x\in[0,1],\,s\in[t,\tau^v(t;x)]\right\},  \label{eq: D}\\
  \D^{\prime}(t) &= \left\{(x,s):\,x\in[0,1],\,s\in[t,\tau^v(t;x))\right\}.  \label{eq: D'}
\end{align}
The state dependences of $\lambda^u$ and $\lambda^v$ in (\ref{tau^u})-(\ref{tau^v}) are included so that the same definitions can be re-used in Section \ref{sec: quasilinear}. 
The following lemma formalizes the above remark that the solution is independent of the control input for a certain amount of time.
\begin{lemma} \label{lemma semilinear determinate set}
Under Assumption \ref{assumption semilinear}, the Cauchy problem consisting of (\ref{w_t})-(\ref{uBC}) with initial conditions (\ref{uvIC})-(\ref{XIC}), i.e., without boundary condition (\ref{vBC}), has a unique bounded solution $u(x,t)$ for $(x,t)\in\D(0)$,   $v(x,t)$ for $(x,t)\in\D^{\prime}(0)$ and $X(t)$ for $t\in [0,\tau^v(0;0)]$, that continuously depends on $(w_0,X_0)$. That is, the solution on these sets is independent of the control input $U(t)$ for $t\geq 0$.
\end{lemma}
\begin{pf}
The proof is  similar to the proof of the first part of Theorem 2 in \cite{strecker2017output}. The difference is that here we have an ODE at the uncontrolled boundary compared to a static boundary condition in \cite{strecker2017output}, which does not affect the determinate sets themselves.  See also  \cite{li2000semi,li2003exact,bressan2000hyperbolic} for a more general discussion of determinate sets (sometimes also called domains of determinacy).
\end{pf}
The solution of $u(x,t)$ for $(x,t)\in\D(0)$ in Lemma \ref{lemma semilinear determinate set} includes the state on the characteristic line $u(\cdot;\tau^v(0;\cdot))$. After shifting time, this implies that the future states $u(\cdot,\tau^v(t;\cdot))$ and $X(\tau^v(t;0))$ are predictable based on the state $w(\cdot,t)$ alone, independently of the  input $U(t)$.

As is usual for hyperbolic PDEs, the state satisfies an ODE along the characteristic lines:
\begin{lemma}\label{lemma semilinear v_x}
For given $t$, the state on the characteristic lines $(x,\tau^v(t;x))$, $x\in[0,1]$, satisfies
\begin{align}
\frac{d}{dx} v(x,\tau^v(t;x)) &= -\frac{f^v(x,(u,v)(x,\tau^v(t;x)))}{\lambda^v(x)}, \label{vtau_x}\\
v(x,\tau^v(t;1)) &= U(t).
\end{align}
\end{lemma}
\begin{pf}
See  \cite[Theorem 2]{strecker2017output}.
\end{pf}
Importantly, Equation (\ref{vtau_x}) can be solved in either $x$-direction. This will be exploited  in the following subsection, where the desired boundary value for $v$ at $x=0$ is designed as a virtual input, and a corresponding actual actual input $U(t)$  is constructed.

\subsection{Sufficient condition for global asymptotic stability}
In order to see that system (\ref{w_t})-(\ref{XIC}) is in fact easier to control via the boundary value of $v(0,\cdot)$, it is insightful to solve (\ref{w_t}) for $w_x$, swapping the roles of $t$ and $x$ and introducing the new input $U^*(t)$. This gives the system
\begin{align}
w_x(x,t)& = \left(\Lambda(x,w(x,t))\right)^{-1}\left(w_t(x,t)-F(x,w(x,t)) \right), \label{w_x}\\
\dot{X}(t) &= f^0(X(t),U^*(t),t), \\
u(0,t) &= g^0(X(t),U^*(t),t),\\
v(0,t) & = U^*(t), \label{vBC2}\\
w(\cdot,0)&=w_0=\left( \begin{matrix}u_0&v_0 \end{matrix} \right)^T,  \\
X(0) &= X_0. \label{XIC2} 
\end{align}
The new input $U^*$ now enters at the common inflow boundary condition ($x=0$) for both distributed states $u$ and $v$, and as an input for the ODE that in turn also affects only the boundary condition for $u$ at $x=0$.

Setting $U^*(t)=K(X(t),t)$ solves the problem of asymptotically stabilizing system (\ref{w_x})-(\ref{XIC2}). 
\begin{lemma}\label{lemma convergence semilinear}
Under Assumption \ref{assumption semilinear}, if $U^*(t)=K(X(t),t)$ for all $t\geq \tau^v(0;0)$, then (\ref{w_x})-(\ref{XIC2}) has a globally asymptotically stable equilibrium at the origin.
\end{lemma}
\begin{pf}
Because of Lipschitz-continuity of the data and the fact that $\Lambda$ is independent of the state, the solution cannot blow up in finite time, i.e., the solutions exists for all $(x,t)\in[0,1]\times[0,\infty)$.
Similar to \cite[Theorem 6]{strecker2021quasilinear2x2}, one can show that there exists some constant $\gamma\geq 1$ such that for all $T\geq 0$,
\begin{equation}
\sup_{x\in[0,1],\,t\geq \tau^u(T;x)} \|w(x,t)\|_{\infty} \leq \gamma\,\sup_{t\geq T}\|w(0,t)\|_{\infty}. \label{condition convergence}
\end{equation}
 For $t< \tau^v(0;0)$, Lemma \ref{lemma semilinear determinate set} ensures that both $w(0,t)$ and $X(t)$ remain bounded and the bounds depend continuously on the initial condition.

For all $t\geq \tau^v(0;0)$, the choice of  $U^*$ ensures that (\ref{closed loop stable semilinear}) is satisfied. Thus, $X(t)$ remains bounded, and the bound during transients can be made arbitrarily small by making the initial condition small. Lipschitz-continuity of $g^0$ and $K$ combined with conditions (\ref{f0(0)=0}) and (\ref{K(0)=0}), imply $\|w(0,\cdot)\|_{\infty}\leq \gamma^{\prime} \|X(\cdot)\|_{\infty}$ for some $\gamma^{\prime}>0$, which implies stability of the whole system. Similarly,  convergence of $X$ to the origin due to asymptotic stability of (\ref{closed loop stable semilinear}), implies that for all $\epsilon>0$ there exists a $T\geq 0$ such that $\|X(t)\|_{\infty}\leq \epsilon $ for all $t\geq T$. Hence $\sup_{t\geq T} \|w(0,t)\|\leq \gamma^{\prime} \epsilon$. Due to (\ref{condition convergence}), this implies that $w$ convergences to the origin as well.
\end{pf}

\subsection{Control law} \label{sec: controller semilinear}
\begin{figure}[htbp!]\centering
\includegraphics[width=.8\columnwidth]{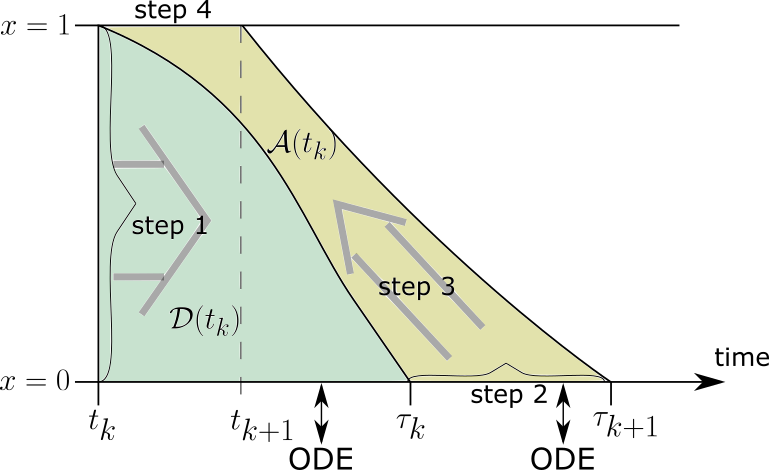}
\caption{The steps of evaluating the controller at time $t_k$, at the example of the sampled-time controller for quasilinear systems:
1) predict state $w$ on determinate set $\D(t_k)$ (shaded in green) and $X$ up to time $\tau_k$; 2) design virtual input $U^*$ 3) solve target system $(\bar{w}^*,\bar{X}^*)$ \emph{backwards} over domain $\A(t_k)$ (yellow); 4) set input to $U(t)=\bar{v}^*(1,t)$. In the semilinear case with  continuous time control law, the domain $\A(t_k)$ reduces to the line $(x,\tau^v(t_k;x))$, $x\in[0,1]$. }
\label{fig: controller steps}
\end{figure}
With the preparations above, we are  in position to synthesize the controller. The  main steps for evaluating the controller are also summarized in Figure \ref{fig: controller steps} and Algorithm \ref{algorithm_semilinear} below.
\subsubsection*{Step 1: Predict state}
 Lemma \ref{lemma semilinear determinate set} can be exploited to establish a prediction of the state on the determinate set $\D(t)$. To this end, denoting the prediction by $(\bar{w},\,X)$ in order to distinguish it from the actual state, we simply copy the dynamics (\ref{w_t})-(\ref{uBC}) but use the current state, $(w(\cdot,t),X(t))$, as the initial condition:
\begin{align}
\bar{w}_t(x,s) &= \Lambda(x,\bar{w}(x,s))\,\bar{w}_x(x,s) + F(x,\bar{w}(x,s)), \label{wbar_t}\\
\dot{\bar{X}}(s) &= f^0(\bar{X}(s),\bar{v}(0,s),s), \label{Xbar_t}\\
\bar{u}(0,s) &= g^0(\bar{X}(s),\bar{v}(0,s),s), \label{ubarBC}\\
\bar{w}(\cdot,t)&= w(\cdot,t) \label{wbarIC}\\
\bar{X}(t) &= X(t). \label{XbarIC}
\end{align}
As stated by Lemma \ref{lemma semilinear determinate set}, we can solve (\ref{wbar_t})-(\ref{XbarIC}) over the determinate set $\D(t)$ to obtain the predictions $\bar{u}(\cdot,\tau^v(t;\cdot))$ and $\bar{X}(\tau^v(t;0))$. Assuming no model uncertainty, these predictions are equal to the values that the actual state will attain. For comparison, estimates on prediction errors when there is uncertainty in parameters and measurements are given in \cite{strecker2021quasilinear2x2}, although only for static boundary conditions.

\subsubsection*{Step 2: Design virtual input $U^*$}
By Lemma \ref{lemma convergence semilinear}, the control objective of asymptotically stabilizing the system is achieved if $v(0,t)$ becomes equal to $K(X(t),t)$. The earliest time at which $U(t)$ has an effect on $v(0,\cdot)$ is given by $\tau^v(t;0)$, and this is also the latest time up to which we can predict $\bar{X}$. Therefore, at each time $t$ we set the future value of the virtual input to
\begin{equation}
U^*(\tau^v(t;0)) = K(\bar{X}(\tau^v(t;0)),\tau^v(t;0)).  \label{U* semilinear}
\end{equation}

\subsubsection*{Step 3: Construct input $U(t)$}
As stated in Lemma \ref{lemma semilinear v_x}, $v(x,\tau^v(t;x))$ satisfies an ODE in $x$ that can in principle be solved in either $x$-direction. The target boundary value that $v(0,\tau^v(t;0))$ shall satisfy is given by (\ref{U* semilinear}). Thus, we can compute the required control input $U(t)$ by solving  the target dynamics
\begin{align}
\frac{d}{dx} \bar{v}^*(x,\tau^v(t;x)) &= -\frac{f^v(x,(\bar{u},\bar{v}^*)(x,\tau^v(t;x)))}{\lambda^v(x)}, \label{vtau_x2}\\
\bar{v}^*(0,\tau^v(t;0)) &= U^*(\tau^v(t;0)),  \label{vtau BC 2}
\end{align}
where we use the prediction of $\bar{u}(x,\tau^v(t;x))$ from Step 1, over the domain $x\in[0,1]$, and set
\begin{equation}
U(t)=\bar{v}^*(1,\tau^v(t;1)). \label{U semilinear}
\end{equation}
 Assuming exact predictions, we have that $v(0,\tau^v(t;0))=U^*(\tau^v(t;0))$ if and only if $U(t)=\bar{v}^*(1,\tau^v(t;1))$  (due to uniqueness of the solution to (\ref{vtau_x2}); see also the proof of Theorem 3 in \cite{strecker2017output}).

\subsubsection*{Control algorithm} 
The steps for evaluating the state feedback controller for semilinear systems at each time $t$ are summarized in the following algorithm.
\begin{algorithm}[H]
\begin{algorithmic}[1]
\REQUIRE state $(w(\cdot,t),X(t))$ \\
\ENSURE control input $U(t)$ \\ ~

\STATE predict the states $\bar{u}(\cdot,\tau^v(t;\cdot))$ and $\bar{X}(\tau^v(t;0))$ by solving (\ref{wbar_t})-(\ref{XbarIC}) over the determinate set $\D(t)$
\STATE set $U^{*}$ as in (\ref{U* semilinear}).
\STATE solve  target dynamics (\ref{vtau_x2})-(\ref{vtau BC 2})  and set $U(t)$ as in (\ref{U semilinear})
\end{algorithmic}
\caption{Control algorithm for semilinear systems }
\label{algorithm_semilinear}
\end{algorithm}

\begin{theorem}
Under Assumption \ref{assumption semilinear}, the system consisting of (\ref{w_t})-(\ref{XIC}) in closed loop with $U(t)$ as constructed by Algorithm \ref{algorithm_semilinear} has a globally asymptotically stable equilibrium at the origin.
\end{theorem}
\begin{pf}
The construction in Algorithm \ref{algorithm_semilinear} ensures that $v(0,t)=K(X(t),t)$ for all $t\geq \tau^v(0;0)$. Therefore, Lemma \ref{lemma convergence semilinear} implies asymptotic stability. 
\end{pf}

\begin{remark}
The assumption of global Lipschitz continuity or global boundedness of $\Lambda^{-1}$ can be restrictive in some cases. If these assumptions hold only locally, similar local results for small initial conditions can be obtained. Briefly, by making the initial condition small one can ensure that the solution remains in a ball around the origin on which a uniform Lipschitz-condition holds. In certain cases, it is also possible to obtain similar results when the assumptions hold only on some not necessarily bounded subsets of the state space. See for instance \cite{gugat2003global}, where the control law keeps the states of the Saint-Venant equations in the  so-called subcritical range.
\end{remark}

\section{Quasilinear systems} \label{sec: quasilinear}
The control design for quasilinear systems is based on the same fundamental ideas as the one for  semilinear systems. However, the need for at least Lipschitz-continuity of the solution requires several changes. Among others, the controller from  Section \ref{sec: controller semilinear} can lead to discontinuities, in which case the solution of a quaslinear system ceases to exist. To avert this, the virtual input must be steered towards the target where it becomes equal to the controller for the nonlinear ODE in a ``slow'' and continuous fashion. Moreover, as will be made precise in Section \ref{sec: existence quasilinear},  it must be ensured that $\|\dot{X}\|_{\infty}$ remains sufficiently small. While there might exist more sophisticated conditions which ensure  that this is the case, it is most straightforward to ensure if the initial conditions are kept small.

Since the compatibility conditions for Lipschitz-continuous solutions imply that $U(t)$ must be equal to $\lim_{x\rightarrow 1} v(x,t)$, the control input is always uniquely determined by the state at that time. Therefore, we propose a sampled-time control law with sampling period $\theta$, where at each sampling instance $t_k=k\,\theta$, $k\in\N$, the control inputs
are pre-computed for the whole time interval $[t_k, t_{k+1}]$. This builds on the work in \cite{strecker2021quasilinear2x2}, where there is no ODE at uncontrolled boundary.

\subsection{Assumptions}
For the quasilinear systems considered in this section, we make the following stronger assumptions.
\begin{assumption}\label{assumption quasilinear}
The speeds are bounded from below and are Lipschitz in the state,
\begin{align}
\sup_{x\in[0,1], w\in\R^2} \min\{ (\lambda^u(x,w))^{-1},\,(\lambda^v(x,w))^{-1} \} &\leq l_{\Lambda^{-1}}, \\
\operatorname*{ess\,sup~max}_{x\in[0,1],w\in\R^2\quad\quad} \{ \|\partial_u \Lambda(x, w) \|_{\infty},\|\partial_v \Lambda(x, w)\|_{\infty} \}&\leq l_{\Lambda} .  \label{bound Lambda_w}
\end{align}
 The nonlinear functions $F$, $f^0$, $g^0$  satisfy the same  Lipschitz-conditions as in (\ref{boundF_w})-(\ref{f0(0)=0}) and there exists a controller satisfying (\ref{bound K_X})-(\ref{K(0)=0}) such that the closed-loop system (\ref{closed loop stable semilinear}) has a  locally asymptotically stable equilibrium at the origin. The initial condition is Lipschitz-continuous and satisfies the compatibility condition
 \begin{equation}
 u_0(0) = g^0(X_0,v_0(0),0).
 \end{equation}
\end{assumption}

\subsection{Existence and predictability of solution} \label{sec: existence quasilinear}
Similar to Lemma \ref{lemma semilinear determinate set}, the solution to (\ref{w_t})-(\ref{XIC}) is  predictable on the determinate sets also in the quasilinear case, although we  need to additionally assume that the initial conditions are small here in order to exclude blow-up of the solution.
\begin{lemma} \label{lemma quasilinear determinate set}
Under Assumption \ref{assumption quasilinear}, there exist $\epsilon>0$ and $\epsilon^{\prime}>0$ such that if $\max\{\|w_0\|_{\infty},\,\|X_0\|_{\infty}\}\leq \epsilon$ and $\max\{\|\partial_t f^0\|_{\infty},\|\partial_t g^0\|_{\infty}, \|w_t(\cdot,0)\|_{\infty},\,\|\dot{X}(0)\|_{\infty}\}\leq \epsilon^{\prime}$, then  the Cauchy problem consisting of (\ref{w_t})-(\ref{uBC}) with initial conditions (\ref{uvIC})-(\ref{XIC}), has a unique bounded solution $w(x,t)$ for $(x,t)\in\D(0)$ and $X(t)$ for $t\in [0,\tau^v(0;0)]$.
\end{lemma}
\begin{pf}
The proof is overall similar to the proof of Theorem 5 in \cite{strecker2021quasilinear2x2}. Briefly, in the quasilinear case one also needs to keep the time derivatives of the PDE states, $w_t$, small. Since these satisfy a PDE that is quadratic in the state, they can blow up in finite time even if the state $w$ remains bounded. It is possible to maintain small $\|w_t\|_{\infty}$ if $\|w_t(0,t)\|_{\infty}$, i.e., the time-derivative at the boundary $x=0$,  is kept sufficiently small. This is more complicated here due to the coupling to the ODE. 

 By virtue of (\ref{vBC}), we have
\begin{equation}\begin{aligned}
\partial_t u(0,t) = &\phantom{ + } \partial_X g^0(X(t),v(0,t),t)\times\dot{X}(t) \\
& + \partial_v g^0(X(t),v(0,t),t) \times\partial_t v(0,t) \\
&+ \partial_t g^0(X(t),v(0,t),t). \label{u_t(0,t)}
\end{aligned}\end{equation}
Using techniques like those  in \cite[proof of Theorem 5]{strecker2021quasilinear2x2} or \cite[Appendix A]{strecker2017output}, one can show that $\max\{\sup_{(x,t)\in\D(0)} \|w(x,t)\|_{\infty}, \sup_{t\in[0,\tau^v(0;0)]}\|X(t)\|_{\infty}\} \leq \gamma \times \max\{ \|w_0\|_{\infty}, \|X_0\|_{\infty}\}$ for some $\gamma\geq 1$. Thus, we can make $\|\dot{X}\|_{\infty}$, as given by the norm of the right-hand side of (\ref{X_t}), arbitrarily small by making $\epsilon$ and $\|\partial_t f^0\|_{\infty}$ small. Similar to \cite{strecker2021quasilinear2x2}, one can also keep $|\partial_t v(0,t)|$ arbitrarily small by making $\epsilon^{\prime}$ and $|\partial_t u(0,s)|$ for all $s<t$  small. In summary, we can make all terms on the right-hand side of (\ref{u_t(0,t)}) sufficiently small by making $\epsilon$ and $\epsilon^{\prime}$ small, so that $\|w_t(0,t)\|_{\infty}$ is small, excluding blow-up of $\|w_t\|_{\infty}$.

In contrast to Lemma \ref{lemma semilinear determinate set}, the control input $U(0)$ is uniquely determined by $U(0)=\lim_{x\rightarrow 1}v_0(x)$. Therefore, $v$ can also be predicted over the closed set $\D(0)$. 
\end{pf}

Similar to Lemma \ref{lemma convergence semilinear} in the semilinear case, a sufficient condition for  asymptotic stability as well as existence of the solution for all times can be formulated by use of the virtual input $U^*$. 
\begin{lemma} \label{lemma condition quasilinear}
Suppose Assumption \ref{assumption quasilinear} holds and that $U^*(t)=K(X(t),t)$ for all $t\geq T$ within finite $T\geq 0$. There exist constants $\epsilon>0$ and $\epsilon^{\prime}>0$ such that if $\max\{\|U^*\|_{\infty},\, \|w_0\|_{\infty},\,\|X_0\|_{\infty}\}\leq \epsilon$ and $\max\{\|\partial_t U^*\|_{\infty},\,\|\partial_t f^0\|_{\infty},\|\partial_t K\|_{\infty},\,\|\partial_t g^0\|_{\infty}, \,\|w_t(\cdot,0)\|_{\infty},$ $\|\dot{X}(0)\|_{\infty}\}\leq \epsilon^{\prime}$, then the solution to (\ref{w_x})-(\ref{XIC2}) exists for all $(x,t)\in[0,1]\times[0,\infty)$ and is bounded by $\gamma\epsilon$ for some $\gamma\geq 1$. Moreover, the solution asymptotically converges to the origin.
\end{lemma}
\begin{pf}
For times $t\leq T$, the norms of $w(0,t)$ and $X(t)$ can be made small with sufficiently small $\epsilon$. Therefore, and using $\|\partial_t f^0\|_{\infty}\leq \epsilon^{\prime}$,  $\|\dot{X}(t)\|$ remains small. Using  (\ref{u_t(0,t)}) and  $|\partial_t v(0,t)|=|\partial_tU^*(t)|\leq \epsilon^{\prime}$, which holds by assumption,  $\|w_t(0,\cdot)\|_{\infty}$ can be made arbitrarily small as well.

 For times $t> T$, $U^*(t)=K(X(t),t)$ implies that $X$ does not grow beyond a bound that can be made arbitrarily small by keeping $X(T)$ small. Moreover, $X$ converges to the origin. Since $\|w(0,t)\|_{\infty}\leq \max\{l_K,l_{g^0}\}\times \|X(t)\|_{\infty}$, $\lim_{t\rightarrow \infty} \|w(0,t)\|{\infty}=0$ and the whole PDE state $w$ converges to the origin as in Lemma \ref{lemma convergence semilinear}. For the time derivatives, we have $\|\dot{X}(t)\|_{\infty}\leq (l_{f^0_X}+l_{f^0_v}l_K)\|X\|_{\infty} + \|\partial_t f^0\|_{\infty} $. Thus, $\|\dot{X}(t)\|_{\infty}$ can be kept below any given bound due to smallness and convergence of $\|X(t)\|_{\infty}$ and smallness of $\|\partial_t f^0\|_{\infty}$. Using Lipschitz-continuity of $K$, we get $|\partial_t v(0,t)|=|\partial_t U^*(t)| \leq l_K \|\dot{X}(t)\|_{\infty} +  \|\partial_t K\|_{\infty}$ and $ \|\partial_t u(0,t)\|_{\infty}$ can again be bounded via (\ref{u_t(0,t)}). 
 In summary, $\|w(0,\cdot)\|_{\infty}$ can be kept  bounded and $\|w_t(0,\cdot)\|_{\infty}$ can be kept arbitrarily small, which ensures existence of the solution $w$ globally on $[0,1]\times[0,\infty)$  like in \cite[Theorem 4]{strecker2021quasilinear2x2}.
\end{pf}

\subsection{Control law}
We can now modify the steps from Section \ref{sec: controller semilinear} to obtain the controller for quasilinear systems. At each time step $t_k$, the sampled-time controller maps the state $(w(\cdot,t_k),X(t_k))$ into the inputs over the interval $[t_k,t_{k+1}]$.
\subsubsection*{Step 1:  Predict state}
The state on the determinate set $\D(t_k)$ can be predicted as in the semilinear case, i.e., by solving (\ref{wbar_t})-(\ref{XbarIC}) over $\D(t_k)$. The  differences are that the state measurement is sampled at discrete times $t=t_k$, that $\bar{v}$ can also be predicted over the closed domain $\D(t_k)$ because of compatibility condition $U(t_k)=\lim_{x\rightarrow 1}v(x,t_k)$, and that the domain $\D(t_k)$ itself depends on the state $(w(\cdot,t_k),X(t_k))$.

\subsubsection*{Step 2: Design virtual input $U^*$}
The goal is still that $U^*(t)$ becomes equal to $K(X(t),t)$. However, the compatibility conditions $U(t_k)=\lim_{x\rightarrow 1}v(x,t_k)$ and $U^*(\tau_k)=\lim_{t\rightarrow \tau_k} v(0,t)$, where $\tau_k=\tau^v(t_k;0)$, mean that this target value cannot be attained immediately. Moreover, Lemma \ref{lemma condition quasilinear} requires a limit on $\partial_t U^*$. One design for $U^*$ that satisfy these two constraints and also converges to the target value $U^*(t)=K(\bar{X}^*(t),t)$, while  remaining bounded during transients, is
\begin{equation}
U^*(t) = \begin{cases} \bar{v}^k + \delta \times \operatorname{sign}(e^k) \times (t-\tau_k), & t\in[\tau_k , \tilde{t}^k] \\
 K(\bar{X}^*(\tau^v(t;0)),\tau^v(t;0)), & t \geq \tilde{t}^k \end{cases},  \label{U* quasilinear}
\end{equation}
where $\bar{v}^k=\bar{v}(0,\tau_k)$, $e^k=K(\bar{X}(\tau_k)),\tau_k)-\bar{v}^k$ is the tracking error at time $\tau_k$, $\delta\leq\epsilon^{\prime}$ is the desired convergence speed and the time $\tilde{t}^k$ is implicitly defined as where the lines $\bar{v}^k + \delta \times \operatorname{sign}(e^k) \times (t-\tau_k)$ and $K(\bar{X}^*(\tau^v(t;0)),\tau^v(t;0))$, $t\geq \tau_k$, intersect for the first time. Note that the term $K(\bar{X}^*(t),t)$  in (\ref{U* quasilinear})  needs to be  evaluated online at the target state $\bar{X}^*$ when solving the target system (\ref{target 1})-(\ref{target last}) below.

\subsubsection*{Step 3: Construct inputs $U(t)$}
Again similar to the semilinear case, the control inputs are constructed by starting with the virtual input $U^*$ and solving the target dynamics along the characteristic lines backwards, relative to the direction the original input propagates. In the quasilinear case, this involves more than solving an ODE as in (\ref{vtau_x2}). 
To shorten notation, we align the time-argument of the target state $(\bar{w}^*,\bar{X}^*)$ with $\tau^v(t;x)$, i.e., $\bar{w}^*(x,t)$ is the target  for the future state $w(x,\tau^v(t;x))$ and $\bar{X}^*(t)$ is the target for $X(\tau^v(t;x))$.
The inputs can be constructed by solving the following system over the domain $(x,t)\in[0,1]\times [t_k,t_{k+1}]$:
\begin{align}
\partial_t \bar{u}^{*}(x,t) =&\ \partial_t \tau^v(t;x)\times \bar{u}^{*}_{\partial t}(x,t), \label{target 1}\\
 \partial_x \bar{v}^{*}(x,t) =& -\frac{f^v(x,\bar{w}^{*}(x,t))}{\lambda^v(x,\bar{w}^{*}(x,t))} , \label{target 2}\\
\partial_t \bar{u}^{*}_{\partial t}(x,t) =& -\mu \times\partial_x\bar{u}^{*}_{\partial t}(x,t)  + \nu \times \left[c_1 \, (\bar{u}^{*}_{\partial t})^2 \right. \nonumber \\
&\left. ~~~~+ c_2\,\bar{u}^{*}_{\partial t}\,\bar{v}^{*}_{\partial t}+c_3 \,\bar{u}^{*}_{\partial t}+c_4\,\bar{v}^{*}_{\partial t}\right],\label{target 3}\\
\partial_x \bar{v}^{*}_{\partial t}(x,t)= & -\frac{1}{\lambda^v} \times \left[c_5\, \bar{u}^{*}_{\partial t}\,\bar{v}^{*}_{\partial t} + c_6\,(\bar{v}^{*}_{\partial t})^2 \right.\nonumber \\
& ~~~~~~~~~~~ \left. + c_7\,\bar{u}^{*}_{\partial t} + c_8\,\bar{v}^{*}_{\partial t} \right] \label{target 4},\\
\partial_t \tau^v(t;x) =&\ 1 - \int_x^1\frac{\partial_u \lambda^v\times \bar{u}^{*}_{\partial t} + \partial_v \lambda^v\times \bar{v}^{*}_{\partial t}}{\left(\lambda^v\right)^2} d\xi, \\
\partial_t \bar{X}^*(t) =&\ \partial_t \tau^v(t;0)\times \bar{X}^*_{\partial t}(t),  \label{target 5}
\end{align}
with boundary conditions
\begin{align}
\bar{v}^{*}(0,t) =&\ \phantom{ \partial_t } U^{*}(\tau^v(t;0)),   \label{target vBC}\\
\bar{v}^{*}_{\partial t}(0,t) =&\ \partial_t U^{*}(\tau^v(t;0)),  \label{target v_tBC} \\
\bar{u}^{*}_{\partial t}(0,t) =&\  \partial_X g^0\times {\bar{X}}^*_{\partial t}(\tau^v(t;0)) \nonumber \\
& + \partial_v g^0\times  \partial_t U^{*}(\tau^v(t;0)) + \partial_t g^0, \label{target uBC}\\ 
{\bar{X}}^*_{\partial t}(t) =&\ f^0(\bar{X}^*(t),  U^{*}(\tau^v(t;0)), \tau^v(t;0)),  \label{target X}
\end{align}
and initial conditions
\begin{align}
\bar{X}^*(t_k) =&\ \bar{X}(\tau_k),\\
\bar{u}^{*}(\cdot,t_k) =&\ \bar{u}(\cdot,\tau_k),\\
\bar{u}_{\partial t}^{*}(\cdot,t_k) =&\ \bar{u}_t^k(\cdot,\tau_k). \label{target last}
\end{align}
Here, the functions $\mu$, $\nu$ and $c_1$ through $c_8$, as well as more details on the derivation behind this system,  are given in \cite[Section III.D]{strecker2021quasilinear2x2}, and  we omitted the arguments of several functions for readability. The subscripts notation $\bar{w}^*_{\partial t}$ is used for the target of the time derivative $w_t(x,\tau^v(t;x))$, where in general $\partial_t\bar{w}^*(x,t) \neq \bar{w}^*_{\partial t}(x,t)$ since $\tau^v$ is time-varying (see also (\ref{target 1}) and (\ref{target 5})). Note that in (\ref{target v_tBC})-(\ref{target uBC}), $\partial_t U^{*,k}(\tau^v(t;0))$ denotes the partial derivative of $U^{*,k} $ with respect to time evaluated at time $\tau^v(t;0)$, not the total derivative of $U^{*,k}(\tau^v(t;0))$ with respect to $t$. 

Because of uniqueness of the solution of (\ref{target 1})-(\ref{target last}), the control inputs such  that $v(0,\tau^v(t;0))$ becomes equal to $U^*(\tau^v(t;0))$ for all $t\in[t_k,t_{k+1}]$ are
\begin{equation}
U(t) = \bar{v}^{*}(1,t). \label{U quasilinear}
\end{equation}

\subsubsection*{Control algorithm} 
The steps for evaluating the state feedback controller for quasilinear systems at each time step $t_k$ are summarized in the following algorithm.
\begin{algorithm}[H]
\begin{algorithmic}[1]
\REQUIRE state $(w(\cdot,t_k),\,X(t_k))$ \\
\ENSURE control input $U(t)$ for $t\in[t_k,t_{k+1}]$ \\ ~

\STATE predict the states $\bar{w}(\cdot,\tau^v(t_k;\cdot))$ and $\bar{X}(\tau_k)$ by solving (\ref{wbar_t})-(\ref{XbarIC}) over the determinate set $\D(t_k)$
\STATE set $U^{*}$ as in (\ref{U* quasilinear}).
\STATE solve  target dynamics (\ref{target 1})-(\ref{target last})  and set $U(t)$ as in (\ref{U quasilinear})
\end{algorithmic}
\caption{Control algorithm for quasilinear systems }
\label{algorithm_quasilinear}
\end{algorithm}

\begin{theorem}
Suppose Assumption \ref{assumption quasilinear} holds. There exists a constant $\delta^{\text{max}}>0$ such that for all $\delta$ satisfying $0<\delta\leq \delta^{\text{max}}$, there exist constants $\epsilon>0$ and $\epsilon^{\prime}>0$ such that if $\max\{\|w_0\|_{\infty},\,\|X_0\|_{\infty}\}\leq \epsilon$ and $\max\{\|\partial_t f^0\|_{\infty},\|\partial_t K\|_{\infty},\,\|\partial_t g^0\|_{\infty}, \,\|w_t(\cdot,0)\|_{\infty},$ $\|\dot{X}(0)\|_{\infty}\}\leq \epsilon^{\prime}$, then the system consisting of (\ref{w_t})-(\ref{XIC}) in closed loop with $U(t)$ as constructed by Algorithm \ref{algorithm_quasilinear} has a locally asymptotically equilibrium at the  origin.
\end{theorem}
\begin{pf}
As in  Lemma \ref{lemma condition quasilinear}, one can make  $\epsilon$ and $\epsilon^{\prime}$ sufficiently small such that $|\frac{d}{dt}K(\bar{X}^*(\tau^v(t;0)),\tau^v(t;0))|\leq 0.5\,\delta$ for sufficiently long, so that the lines $\bar{v}^k + \delta \times \operatorname{sign}(e^k) \times (t-\tau_k)$ and $K(\bar{X}^*(\tau^v(t;0)),\tau^v(t;0))$ intersect within finite time. That is, $U^*$ as given in (\ref{U* quasilinear}) becomes equal to  $K(\bar{X}^*(\tau^v(t;0)),\tau^v(t;0))$ within finite time. This, combined with the  construction in Algorithm \ref{algorithm_quasilinear}, ensures that $v(0,t)=K(X(t),t)$ within finite time.  Moreover, the design ensures that for all times where $v(0,t)\neq K(X(t),t)$, $v(0,t)=U^*(t)$ satisfies the conditions of Lemma \ref{lemma condition quasilinear}. Therefore, Lemma \ref{lemma condition quasilinear} implies asymptotic stability and that the solution does not cease to exist and remains below a bound depending on the initial condition.
\end{pf}

\section{State  estimation and output feedback control}\label{sec:observer}
\begin{figure}[htbp!]\centering
\includegraphics[width=.5\columnwidth]{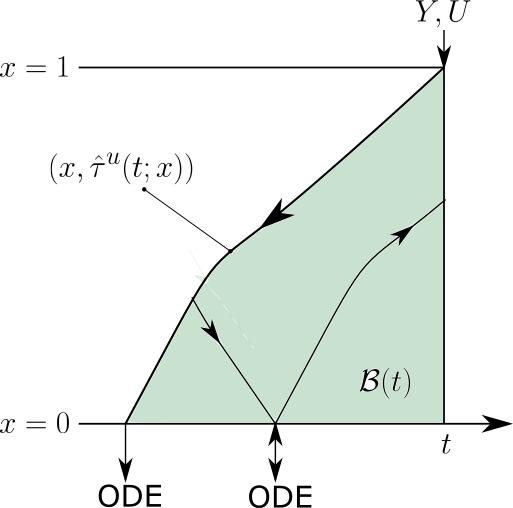}
\caption{Schematic of the observer. First, the state on the characteristic line $(x,\hat{\tau}^u(t;x))$, $x\in[0,1]$ is estimated, starting with the measurement at $x=1$ and \emph{backwards} relative to the propagation direction of $u$. Second, the current state is estimated by solving the  system dynamics over the domain $\B(t)$ (shaded in green).  }
\label{fig:observer}
\end{figure}
Sections \ref{sec: semilinear} and \ref{sec: quasilinear} considered the case where the whole state $(w(\cdot,t),X(t))$ is available to the feedback controller, which is unrealistic in practical applications. Therefore, we present an observer that estimates both the distributed state $w(\cdot,t)$ and the ODE state $X(t)$ from the boundary measurement 
\begin{equation}
Y(t) = u(1,t)  \label{Y}.
\end{equation} 
 Like in \cite{strecker2017output}, the observer design is based on the observation that the measurement $Y(t)$ contains only information about \emph{past} values of the ODE state and the distributed state in the interior of the domain. For instance, the effect of any disturbance affecting $X$ would not be detectable in $Y$ for a certain amount of time. Therefore, we first estimate the past state on the characteristic line along which the measurement evolved, before mapping this estimate into an estimate of the current state via a prediction step. This two-step procedure is outlined in Figure \ref{fig:observer}. We again split the observer design into two cases depending on whether the speeds $\Lambda$ depend on the state or not.

\begin{assumption}\label{assumption observer}
There exists an observer $h$ for the ODE part of the system such that the state estimate $\hat{X}(t)$ satisfying
\begin{align}
\dot{\hat{X}}(t) &= h(\hat{X}(t),u(0,t),v(0,t),t), \\
\hat{X}(0) &= \hat{X}_0, 
\end{align}
converges to the actual value $X(t)$ within  finite time for any initial guess $\hat{X}_0$, input signal $v(0,t)$ and output measurement $u(0,t)$.
\end{assumption}
Some results on finite-time estimation of ODEs are given in \cite{engel2002finitetimeobserver,menard2010finitetimeobserver,lopez2018finitetimeobserver}.
\begin{remark}
The assumption that the observer for the ODE part of the system converges within finite time is conservative and is made only so that  it fits into the framework of exact estimation and exact predictions as in the state feedback case. Algorithmically,  the implementation for the case where the observer for the ODE  converges only asymptotically is exactly the same. However, any remaining estimation error would lead to  prediction error when  mapping the estimates $(\hat{w},\hat{X})$ into the estimate of the current state $(w,X)$ as described in Section \ref{sec: observer second step}. In the case with static boundary conditions, it was possible to show that the closed-loop system is robustly stable in presence of uncertainty in the boundary condition \cite{strecker2021quasilinear2x2}, and one might expect that something similar can be shown for the case considered here. 
\end{remark}

Define the characteristic line along which the measurement $Y(t)$ evolved and the past state on this line  as
\begin{align}
\check{\tau}^u(t;x) &= t - \int_x^1 \frac{1}{\lambda^u(\xi;w(\xi,\check{\tau}^u(t;\xi)))}d\xi, \label{tau_hat}  \\
\check{w}(x,t) &= (\check{u},\check{v})(x,t) =  w(x,\check{\tau}^u(t;x)), \label{wcheck} \\
\check{X}(t) &= X(\check{\tau}^u(t;0)). \label{Xcheck}
\end{align}

\subsection{Estimation of $(\check{w},\check{X})$ - semilinear systems}
If  $\Lambda$ is independent  of the state, one can show that, similar to \cite[Theorem 8]{strecker2017output}, $\check{w}$ satisfies
\begin{align}
\check{u}_x(x,t) &= \frac{f^u(x,\check{w}(x,t))}{\lambda^u(x)}, \label{ucheck_x}\\
\check{v}_t(x,t) &= \hat{\mu}(x) \,\check{v}_x(x,t) + \hat{\nu}(x) \,f^v(x,\check{w}(x,t)) ,\label{vcheck_t}\\
\dot{\check{X}}(t) &= f^0(\check{X}(t),\check{v}(0,t),\check{\tau}^u(t;0)), \label{Xcheck_t}\\
\check{u}(0,t) &= g^0(\check{X}(t),\check{v}(0,t),\check{\tau}^u(t;0)),\label{ucheck0}\\
\check{v}(1,t) &= U(t), \label{vcheck1}
\end{align}
where $\hat{\nu} = \frac{\lambda^u}{\lambda^u+\lambda^v}$ and $\hat{\mu} = \lambda^v \hat{\nu}$. Note that (\ref{ucheck_x}) is an ODE which can again be solved in either $x$-direction. The boundary value of $\check{u}$ at $x=1$ is known, whereas the boundary condition at $x=0$  given by (\ref{ucheck0}) depends on $\check{X}$, which is yet to be estimated. Therefore, we design the observer as a copy (\ref{ucheck_x})-(\ref{vcheck1}) but with the boundary condition (\ref{ucheck0}) replaced by the measurement, and the ODE (\ref{Xcheck_t}) replaced by the observer from Assumption \ref{assumption observer}:
\begin{align}
\hat{u}_x(x,t) &= \frac{f^u(x,\hat{w}(x,t))}{\lambda^u(x)}, \label{uhat_x}\\
\hat{v}_t(x,t) &= \hat{\mu}(x) \, \hat{v}_x(x,t) + \hat{\nu}(x)\, f^v(x,\hat{w}(x,t)) ,\label{vhat_t}\\
\hat{u}(1,t) &= Y(t),\label{uhat1}\\
\check{v}(1,t) &= U(t), \label{vhat1}\\
\hat{w}(\cdot,0) &= \hat{w}_0,\label{what_IC}\\
\dot{\hat{X}}(t) &= h(\hat{X}(t),\hat{u}(0,t),\hat{v}(0,t),\check{\tau}^u(t;0)), \label{Xhat_t}\\
\hat{X}(0) &= \hat{X}_0. \label{Xhat_IC}
\end{align}
Note that observer (\ref{uhat_x})-(\ref{Xhat_IC}) has a cascade structure: The PDE state $\hat{w}$ is estimated via (\ref{uhat_x})-(\ref{what_IC}) without using the boundary condition at $x=0$ or $\hat{X}$ at all; then, the estimates of the unmeasured boundary values  $(\hat{u}(0,t),\hat{v}(0,t))$ are used to estimate the ODE state $\hat{X}$ via (\ref{Xhat_t})-(\ref{Xhat_IC}). That is, the estimate  $\hat{X}$ does not affect the estimate $\hat{w}$ at all.   With this design, the estimates $(\hat{w},\hat{X})$ converge to the actual, past states $(\check{w},\check{X})$ within finite time.
\begin{theorem}\label{thm observer semilinear}
There exists a $T>0$ such that for any initial guess $(\hat{w}_0,\hat{X}_0)$ and any $Y,U$ in $L^{\infty}$, we have $\hat{w}(\cdot,t) = \check{w}(\cdot,t)$ and $\hat{X}(t)= \check{X}(t)$ for all $t\geq T$.
\end{theorem}
\begin{pf}
The PDE-part of the observer, (\ref{uhat_x})-(\ref{what_IC}), converges to $\check{w}$ within finite time exactly  as in the case with static boundary condition in \cite[Theorem 19]{strecker2017output}. That is, within finite time one has exact estimates $(\hat{u}(0,t),\hat{v}(0,t))=(\check{u}(0,t),\check{v}(0,t))$, which are the inputs to the ODE-observer. Therefore,  $\hat{X}$ converges to $\check{X}$ within another finite amount of time due to Assumption \ref{assumption observer}.
\end{pf}

\subsection{Estimation of $(\check{w},\check{X})$ - quasilinear systems}
In quasilinear systems, the speed corresponding to $\hat{\mu}$ in (\ref{vcheck_t}) becomes time varying because the characteristic line $\check{\tau}^u$ depends on the state. Therefore, the observer additionally needs to estimate $\check{\tau}^u$ and its  time-derivative. This is possible by designing the observer as a copy of the dynamics of $\check{w}$ up to the first time-derivatives (in a manner similar to (\ref{target 1})-(\ref{target last})), and again imposing the measurements and its time-derivative  as a boundary value at $x=1$ instead of using the to-be estimated boundary condition at $x=0$. We estimate the PDE state via
\begin{align}
\partial_x \hat{u}(x,t) =& \ \frac{f^u}{\lambda^u}, \label{obs q 1}\\
\partial_t \hat{v}(x,t) =&\ \partial_t \hat{\tau}^u \times\left( \hat{\mu} \times \partial_x \hat{v}(x,t) + \hat{\nu}\times f^v \right)  \label{obs q 2}\\ 
\partial_x \hat{u}_{\partial t}(x,t) =&\  \frac{1}{\lambda^u} \times \left[c_1 \, (\hat{u}_{\partial t})^2+ c_2\,\hat{u}_{\partial t}\,\hat{v}_{\partial t} \right. \nonumber \\
&\left. ~~~~~~~~~+c_3 \,\hat{u}_{\partial t}+c_4\,\hat{v}_{\partial t}\right],\label{obs q 3}\\
\partial_t \hat{v}_{\partial t}(x,t)= &\ \partial_t \hat{\tau}^u \times \left( \hat{\mu}\times\partial_x \hat{v}_{\partial t}(x,t) + \hat{\nu}  \times \left[c_5\, \hat{u}_{\partial t}\,\hat{v}_{\partial t}  \right.\right.\nonumber \\
& ~~~~ \left.\left. + c_6\,(\hat{v}_{\partial t})^2 + c_7\,\hat{u}_{\partial t} + c_8\,\hat{v}_{\partial t} \right]\right) \label{obs q 4},\\
\hat{\tau}^u(t;x) = &\ t - \int_x^1 \frac{1}{\lambda^u(\xi;\hat{w}(\xi,t))}d\xi \\
\partial_t \hat{\tau}^u(t;x) =&\ 1 + \int_x^1\frac{\partial_u \lambda^u\times \hat{u}_{\partial t} + \partial_v \lambda^u\times \hat{v}_{\partial t}}{\left(\lambda^u\right)^2} d\xi, \\
\hat{u}(1,t)=&\ Y(t), \label{obs q uBC}\\
\hat{u}_{\partial_t}(1,t)=&\ \partial_t Y(t), \label{obs q u_tBC}  \\
\hat{v}(1,t)=&\ U(t), \label{obs q vBC}\\
\hat{v}_{\partial_t}(1,t)=&\ \partial_t U(t), \label{obs q v_tBC}
\end{align}
and the ODE state via
\begin{equation}
\dot{\hat{X}}(t)= \partial_t \hat{\tau}^u(t;0)\times h(\hat{X}(t),\hat{u}(0,t),\hat{v}(0,t),\hat{\tau}^u(t;0)),  \label{obs q X}
\end{equation}
with initial guesses
\begin{align}
\hat{w}(\cdot,0)&= \hat{w}_0, & \hat{v}_{\partial t}(\cdot,0) &=  \hat{v}_{\partial t}^0, & \hat{X}(0)&=\hat{X}_0.  \label{obs q IC}
\end{align}

Similar to the semilinear case, this observer converges in finite time. 
\begin{theorem}\label{thm observer quasilinear}
There exist $\epsilon^{\prime}>0$ (which depends on $\max\{\|\hat{w}_0\|_{\infty},\|Y\|_{\infty},\|U\|_{\infty}\}$) and $T>0$ such that if $\| \hat{v}_{\partial t}^0\|_{\infty}\leq \epsilon^{\prime}$, $\|\partial_t Y\|_{\infty}\leq \epsilon^{\prime}$ and $\|\partial_t U\|_{\infty}\leq \epsilon^{\prime}$, then $\hat{w}(\cdot,t) = \check{w}(\cdot,t)$ and $\hat{X}(t)= \check{X}(t)$ for all $t\geq T$.
\end{theorem}
\begin{pf}
Convergence follows exactly as in the semilinear case in Theorem \ref{thm observer semilinear}.  Smallness of $\hat{v}_{\partial t}^0$, $\partial_t Y$ and $\partial_t U$ is sufficient for preventing blow-up of the estimate $\hat{w}$. This can be shown using the concept of semi-global solutions as in  \cite[Theorems 3 and 4]{strecker2021quasilinear2x2}. Briefly, the characteristic lines of (\ref{obs q 1})-(\ref{obs q IC})  all start at either $x=1$ or $t=0$ and all go in the negative $x$-direction. Consequently,  (\ref{obs q 1})-(\ref{obs q IC}) (excluding the part involving $X$)   can be seen as an PDE-ODE system in the negative $x$-direction, and   conditions like in Lemma \ref{lemma condition quasilinear} or  \cite[Theorem 4]{strecker2021quasilinear2x2} can be applied to ensure existence of a solution for $\hat{w}$ on the bounded $x$-horizon $x\in[0,1]$. The estimate of the ODE state, $\hat{X}$, does not blow up and it converges due to Assumption \ref{assumption observer}.
\end{pf}

\subsection{Estimation of current state and output feedback}\label{sec: observer second step}
The observers presented above only estimate the past state $(\check{w},\check{X})$ as defined in (\ref{wcheck})-(\ref{Xcheck}). That estimate of the past state can be mapped to an estimate of the current state, at time $t$, through another prediction step by solving the system dynamics with initial condition given by $(\hat{w},\hat{X})$, i.e.,
\begin{align}
w^{\mathrm{est}}_t(x,s) &= \Lambda(x,w^{\mathrm{est}}(x,s))\,w^{\mathrm{est}}_x(x,s) \nonumber \\
&\qquad  + F(x,w^{\mathrm{est}}(x,s)), \label{west_t}\\
\dot{X}^{\mathrm{est}}(s) &= f^0(X^{\mathrm{est}}(s),v^{\mathrm{est}}(0,s),s), \label{Xest_t}\\
u^{\mathrm{est}}(0,s) &= g^0(X^{\mathrm{est}}(s),v^{\mathrm{est}}(0,s),s), \label{uestBC}\\
w^{\mathrm{est}}(\cdot,\hat{\tau}^u(t;\cdot))&=\hat{w}(\cdot,t),\\
X^{\mathrm{est}}(\hat{\tau}^u(t;0)) &= \hat{X}(t), \label{XestIC}
\end{align}
over the domain 
\begin{align}
 \B(t) &= \left\{(x,s):\,x\in[0,1],\,s\in[\hat{\tau}^u(t;x),t]\right\},  \label{eq: B}
\end{align}
where $\hat{\tau}^u$ can be replaced by $\check{\tau}^u$ in the semilinear case.
Similar to Lemmas \ref{lemma semilinear determinate set} and \ref{lemma quasilinear determinate set}, and assuming in the quasilinear case that $\hat{w}_0$, $\hat{X}_0$, $Y$ and $U$ have sufficiently small infinity-norm, system (\ref{west_t})-(\ref{XestIC}) has a unique solution on $\B(t)$. By virtue of Theorems \ref{thm observer semilinear} and \ref{thm observer quasilinear}, this implies that $(w^{\mathrm{est}}(\cdot,t),X^{\mathrm{est}}(t)) = (w(\cdot,t),X(t))$ for all $t\geq T$ for some sufficiently large, but finite $T$.

This observer can be combined with the state feedback controllers from Sections \ref{sec: semilinear} and \ref{sec: quasilinear}, respectively, to obtain the corresponding output feedback controller that uses only the boundary measurement $Y$ as defined in (\ref{Y}). In the semilinear case, the point $x=1$ can be left out of the domain $\B(t)$ when solving (\ref{west_t})-(\ref{XestIC}) because the control input $U(t)$ is not yet known when the observer is evaluated. Omitting this one point has no influence on the solution of  (\ref{vtau_x2}) when  the control input is constructed. In the quasilinear case, smallness of $U$, $Y$ and their derivatives is also required. This can be guaranteed if $\|w_0\|_{\infty}$, $\|X_0\|_{\infty}$ and $\|\partial_x w_0\|_{\infty}$ are sufficiently small, and $\|\hat{w}_0\|_{\infty}$, $\|\hat{X}_0\|_{\infty}$, $\|\partial_x \hat{w}_0\|_{\infty}$ as well as $\delta$ in  (\ref{U* quasilinear}) are chosen sufficiently small. This can be shown via analysis of the worst-case growth of the time derivatives as performed in \cite{strecker2021quasilinear2x2}.
 In practice it makes sense, depending on the application,  to activate the controller only once the observer has had  time to converge to the actual state.
 
 The observer presented here is somewhat related to the one from \cite{li2008observability} for systems with static boundary conditions, in that it provides exact state estimates within the same time. However, in the approach from \cite{li2008observability}, the estimate of the past state in the interior of the domain is obtained by saving the history of boundary measurements and solving another PDE. Whereas in the approach presented here,  the observers (\ref{uhat_x})-(\ref{Xhat_IC}) or (\ref{obs q 1})-(\ref{obs q IC}), respectively, estimate $(\check{w},\check{X})$ dynamically without the need for saving any measurements.

\section{Simulation example} \label{sec:simulation}
We demonstrate the performance of the proposed controller in numerical simulations of a system with $n=1$ and the parameters
\begin{align}
f^u&= \sin(u+v), & f^v&= \sin(v-u), \\
 \lambda^u &= \begin{cases} 0.5 & x<0.5 \\ x& x\geq 0.5 \end{cases}, &  \lambda^u &= 1+0.5\times(|u|+|v|), \\
 f^0&=X|X|+v(0,t), & g^0&= X+v(0,t), \\
 v_0&= 0.5\times(1+x), & u_0&= -0.5,\\
 X_0&= -1.
\end{align}
The sampling time is set to $\theta=0.5$ and the convergence speed to $\delta=1$. For simplicity, we consider the state-feedback case. As in \cite{strecker2021quasilinear2x2}, the simulation and controller are implemented by first semi-discretizing all PDEs in space using finite differences. The resulting high-order ODEs are then solved in Matlab by use of ode45. For a discretization with 100 spatial elements, which has been used to produce the figures, evaluating the controller takes around $0.1$ seconds on a standard laptop, although the code has not been optimized for performance.

\textbf{Stabilization}~
\begin{figure*}[htbp!]\centering
\includegraphics[width=0.3\textwidth]{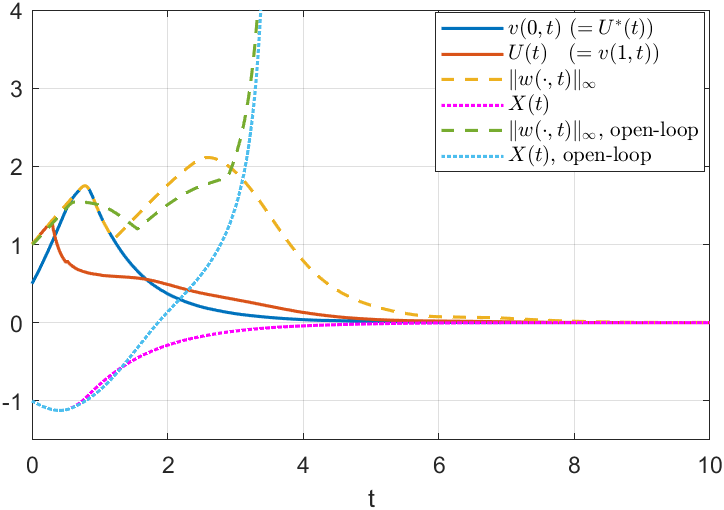}\hfill
\includegraphics[width=0.32\textwidth]{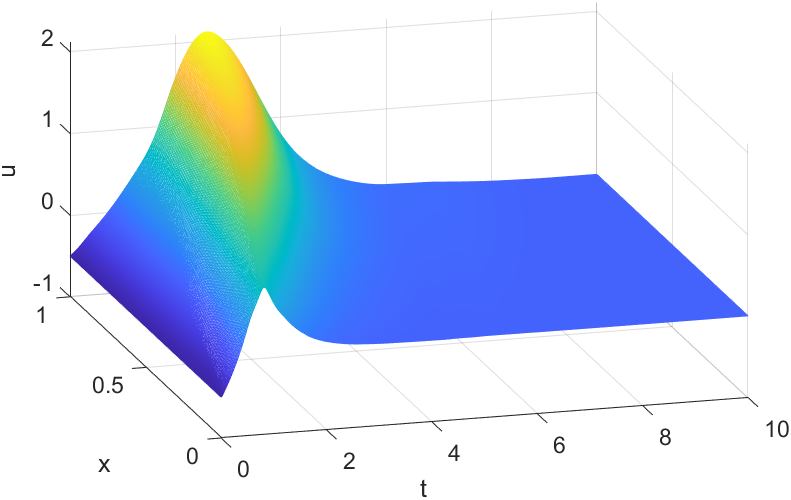}\hfill
\includegraphics[width=0.32\textwidth]{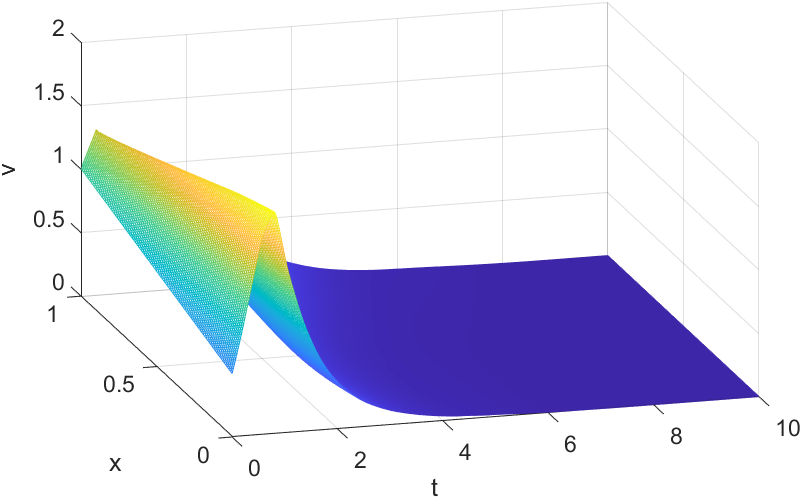}\\
\includegraphics[width=0.3\textwidth]{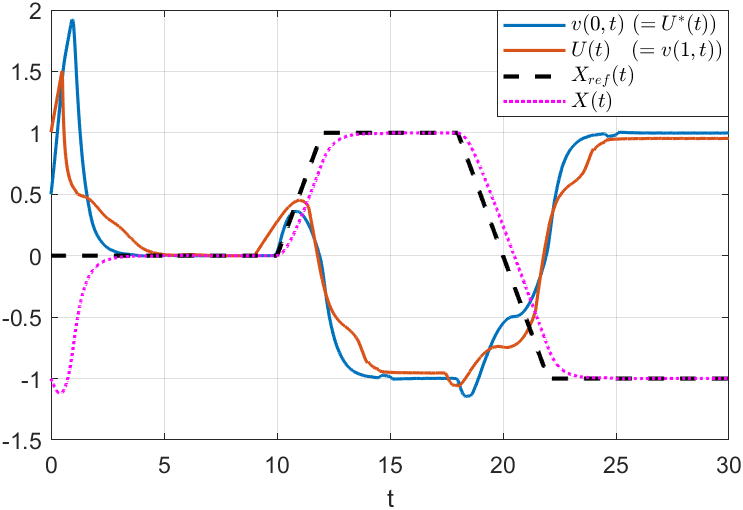}\hfill
\includegraphics[width=0.32\textwidth]{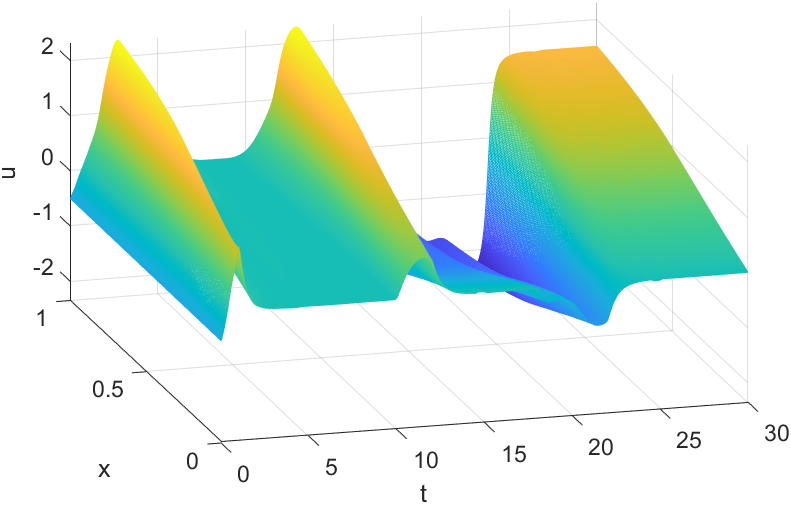}\hfill
\includegraphics[width=0.32\textwidth]{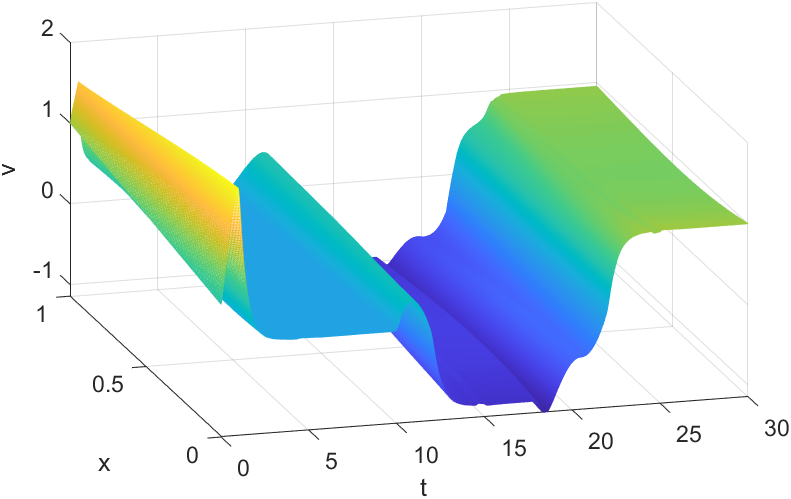}
\caption{Trajectories of the stabilization (top) and tracking (bottom) examples. The top left figure also shows the state norm and trajectory of $X$ for the open-loop case where the control input is held constant at $U(t)=1$. }
\label{fig:examples}
\end{figure*}
Figure \ref{fig:examples} shows the trajectories when the controller 
\begin{equation}
K = -X|X| - X
\end{equation}
is used. Clearly, the term $-X|X|$ is designed to cancel the destabilizing nonlinearity in $f^0$, while the term $-X$ is used to drive the system to the origin. In accordance with theory, after some transients the system asymptotically converges to the origin. For comparison, Figure \ref{fig:examples} also shows the trajectories of $X$ and $\|w(\cdot,t)\|_{\infty}$ for the case where the system is run in open loop with $U(t)=1$ for all times. For this constant input, the solution escapes in finite time at around $t=3.6$.

\textbf{Reference tracking}~
The control design presented in this paper can be modified to achieve  tracking of a reference signal $X_{ref}$. If tracking is the control objective, it is not necessary to assume that the origin is an equilibrium, i.e., assumptions (\ref{F(0)=0})-(\ref{f0(0)=0}) and (\ref{K(0)=0}) can be dropped. For this purpose, the ODE controller is changed to 
\begin{equation} 
K = -X|X| +2\times (X_{ref}(t)-X).
\end{equation}
That is, the controller becomes time-varying. Like before, the term $-X|X| $ is to cancel the destabilizing nonlinearity and the second term $2\times (X_{ref}(t)-X)$ is designed to achieve reference tracking.  The simulated  trajectories are shown in the bottom row of Figure \ref{fig:examples}. The controller achieves good reference tracking and the solution of remains bounded at all times. 

\section{Conclusions} \label{sec:conclusion}
We generalized a recently developed boundary controller for $2\times 2$ quasilinear hyperbolic systems to the case where the same type of PDE is coupled with an ODE at the uncontrolled boundary. Based on predictions over the maximum determinate sets, it is possible to treat the future control input to the ODE subsystem as an virtual input, before the actual inputs are computed by feeding the virtual inputs through the inverted PDE dynamics. A detailed comparison of semilinear and general quasilinear systems is given. In particular, the control inputs can be chosen more freely in the semilinear case and less attention needs to be put into ensuring that the solution does not blow up in finite time.

It is noteworthy that the generalization from $2\times 2$ hyperbolic systems with static boundary conditions to systems coupled with an ODE is relatively straightforward despite the complex  dynamics. Some other variations that have been published in the past for semilinear systems, such as to a linear network structure \cite{strecker2017seriesinterconnections} and to bilateral control \cite{strecker2017twosided}, have also been without major complications once the design for basic $2\times 2$ hyperbolic systems was known. Overall, these results demonstrate the flexibility of the approach. 

We also presented an observer that estimates the distributed PDE state and the unmeasured ODE state from measurements at the actuated boundary only. The observer has a cascade structure, where the PDE state is estimated independently of the estimate of the ODE state, but the estimate of the PDE state is used to estimate the ODE state. This structure can be exploited if the hyperbolic PDE is coupled to a different type of system, such as another PDE. Another situation of interest might be $2\times 2$ hyperbolic PDEs coupled to an uncertain  ODE. If, for instance, an adaptive observer existed for the ODE itself (without the coupling to the PDE), it is straightforward to include such an adaptive observer in the design presented here, although rigorous conditions for existence of the solutions might be challenging. 

\bibliographystyle{plain}        
\bibliography{references}           

\end{document}